\documentclass[11pt]{amsart}
\usepackage[utf8]{inputenc}
\usepackage[T1]{fontenc}
\usepackage{lmodern,amsmath,amssymb,amsthm,mathrsfs,microtype, xcolor,enumerate}
\usepackage[normalem]{ulem}
\usepackage[margin=1in]{geometry}
\usepackage[colorlinks=true,linkcolor=blue,citecolor=blue,urlcolor=blue]{hyperref}
\newcommand{\R}{\mathbb R}

\newcommand{\HC}{\operatorname{HC}}
\newcommand{\supp}{\operatorname{supp}}
\newcommand{\diam}{\operatorname{diam}}
\newcommand{\hrel}{\overset{h}{\sim}}
\newcommand{\hess}{H_{\mathrm{ess}}}
\newcommand{\htop}{h_{\mathrm{top}}}

\newtheorem{theorem}{Theorem}[section]

\newtheorem{lemma}[theorem]{Lemma}
\newtheorem{corollary}[theorem]{Corollary}
\newtheorem{mainthm}{Theorem}

\theoremstyle{definition}
\newtheorem{definition}[theorem]{Definition}
\theoremstyle{remark}
\newtheorem{remark}[theorem]{Remark}
\numberwithin{equation}{section}
\title[MME for factors of surface diffeomorphisms]{Measures of maximal entropy for monotone factors\\ of smooth surface diffeomorphisms}

\makeatletter
\let\@affiliations\@empty
\newcommand{\affiliations}[1]{\gdef\@affiliations{#1}}

\apptocmd{\@setauthors}{%
	\ifx\@affiliations\@empty\else
	\par\vspace{0.8\baselineskip}%
	{\centering\small\@affiliations\par}%
	\fi}
\makeatother

\author[S. Crovisier]{Sylvain Crovisier\textsuperscript{1}}
\author[A. Garc\'ia-Sassi]{Alejo Garc\'ia-Sassi\textsuperscript{2}}
\author[Y. Lima]{Yuri Lima\textsuperscript{3}}
\author[F. Tal]{F\'abio Tal\textsuperscript{4}}

\affiliations{%
	\textsuperscript{1} LMO, CNRS / Université Paris–Saclay, France.\\
	\emph{Email address:} \texttt{sylvain.crovisier@universite-paris-saclay.fr}\\[0.5ex]
	\textsuperscript{2}IMERL, Universidad de la República, Uruguay --- IME, Universidade de São Paulo, Brasil. \\
	\emph{Email address:} \texttt{alejog@fing.edu.uy}\\[0.5ex]
	\textsuperscript{3}IME, Universidade de São Paulo, Brasil.
	\emph{Email address:} \texttt{yurilima@ime.usp.br}\\[0.5ex]
	\textsuperscript{4}IME, Universidade de São Paulo, Brasil.
	\emph{Email address:} \texttt{fabiotal@ime.usp.br}%
}
\date{}

\begin{document}
\begin{abstract}
We prove that a transitive monotone factor of a smooth
surface diffeomorphism has a unique measure of maximal entropy whenever
its topological entropy is positive. Its entropy is the entropy of a
single homoclinic class of the diffeomorphism, and the factor map is a
measure-theoretic isomorphism on every ergodic measure with positive entropy
of that class. The other homoclinic classes are mapped to
periodic orbits. We give an application to the essential factor of a diffeomorphism
isotopic to a pseudo-Anosov map ---possibly relative to a finite invariant
set--- constructed by Garc\'ia-Sassi \& Tal: the essential factor has a unique measure of maximal entropy, which is Bernoulli and has full support.
\end{abstract}
\maketitle

\section{Introduction}\label{sec:intro}

\newcommand\SY[1]{{\color{blue}#1}}
\newcommand{\sy}[1]{\marginpar{\tiny \YL{#1}}}
\newcommand\YL[1]{{\color{red}#1}}
\newcommand{\yl}[1]{\marginpar{\tiny \YL{#1}}}

A surface diffeomorphism can carry several distinct sources of positive topological
entropy. Some reflect the global topology of the dynamics, while others
are confined to regions that disappear under a natural quotient. The
essential factor construction of Garc\'ia-Sassi \& Tal
\cite{garciatal24chaotic,garciatal26essential} separates these
phenomena by collapsing dynamically bounded continua. The resulting
factor retains the large-scale dynamics and has strong topological
properties, although it need not inherit a differentiable structure
compatible with its dynamics. We ask whether this factor has a unique
measure of maximal entropy.

For $C^\infty$ surface diffeomorphisms, Buzzi, Crovisier \& Sarig
\cite{buzzicrovisiersarig22measures} proved finiteness of the ergodic
measures of maximal entropy, and uniqueness under topological
transitivity. Their homoclinic class theory gives a more precise
statement: each homoclinic class with positive topological entropy supports a unique
{probability measure that maximizes the entropy among the invariant probability measures in the class}, it has
full support on the class and a Bernoulli
structure up to a finite period. These results do not apply to the essential factor, which is only a homeomorphism.
Moreover, the original dynamics can have entropy which the quotient removes.
Consequently, the measures of maximal entropy {of the diffeomorphism} need not be the
appropriate measures to project.

In this work, we identify a {\em unique} homoclinic class carrying the topological entropy of the
factor. Connectedness of the fibres is the crucial property: a fibre
through a typical point is a singleton unless the
projected measure is periodic. 
Combining this fact with the
homoclinic class theory yields a result for arbitrary transitive
monotone factors of smooth surface diffeomorphisms. The essential
factor is an application of this more general statement.

Throughout, a surface is connected, closed and smooth. A continuous
surjection is \emph{monotone} when each of its fibres is connected.
{We only consider probability measures.} A measure of maximal entropy is abbreviated as MME. We use \emph{Bernoulli up to a period} to mean that a measure is isomorphic to the
product of a Bernoulli shift and a cyclic permutation on a finite set.

Let $S$ be {an} orientable surface of genus at least one and let
$f$ be a $C^\infty$ diffeomorphism isotopic to a pseudo-Anosov homeomorphism, possibly relative to a finite $f$-invariant set.
For the isotopy class of the identity, this coincides with diffeomorphisms whose ergodic rotation set has nonempty interior, combining the works 
\cite{franks89realizing,llibremackay91rotation,boyland09transitivity} for the torus case, 
and by \cite{guiheneufmiliton23hyperbolic} or \cite{garciaguiheneuflessa23tracking} for the case of genus at least two. 
The {\em essential factor} constructed in \cite{garciatal24chaotic, garciatal26essential} is denoted by
\[
 \pi:S\to S_{\mathrm{ess}} \ \simeq S,\qquad
 \pi\circ f=g\circ\pi.
\]
It collapses the maximal dynamically bounded continua (a continuum is dynamically bounded if
its compact
lifts to the universal cover have uniformly bounded diameters under all
positive and negative iterates). We use the quoted results that
$S_{\mathrm{ess}}$ is a surface homeomorphic to $S$, $\pi$ is monotone, $g$ is topologically
mixing,
and every nonempty open set of $S_{\mathrm{ess}}$ contains a topological
horseshoe for some iterate of $g$. The geometric ingredients of that
construction include the generalized foliations of
\cite{militon24foliations} and the homotopical stability theory of
\cite{fathi90homotopical}, as well as some forcing theory from \cite{lecalveztal18forcing,lecalveztal22topological}.

For a saddle periodic orbit $\mathcal O$, let $\HC(\mathcal O)$ denote its homoclinic
class. An ergodic measure is hyperbolic if it has one positive and one negative Lyapunov exponent. We recall,
{by the Ruelle inequality,} that any ergodic measure with positive entropy is hyperbolic. We write $\mu\hrel\HC(\mathcal O)$ when the ergodic hyperbolic measure
$\mu$ is homoclinically related to $\mathcal O$ (see Definition \ref{def:homoc_relat_meas}). The period of $\HC(\mathcal O)$ is the greatest
common divisor of the periods of the periodic orbits related to $\mathcal O$.

\begin{mainthm}[\bf {MME for the essential factor}]\label{thm:A}
Let $f:S\to S$ be a $C^\infty$ surface diffeomorphism, isotopic to a pseudo-Anosov homeomorphism relative to a finite set, and let $g$ be its essential factor. Then $g$ has a unique MME, which is Bernoulli and has full support.
More specifically, there is a unique homoclinic class
$\hess$ of $f$ whose image {by $\pi$} is not a finite periodic orbit, and it satisfies:
\begin{enumerate}[\rm (1)]
\item $\pi(\hess)=S_{\mathrm{ess}}$ and
$\htop(g)=\htop(f|_{\hess})
 =\sup\{h_\mu(f):\mu\text{ ergodic and }\mu\hrel\hess\}$. 
\item If $\mu^*$ is the unique MME of $f|_{\hess}$, then $\nu^*=\pi_*\mu^*$ is the unique MME of $g$,
and $\pi$ induces an isomorphism between $(f,\mu^*)$ and $(g,\nu^*)$.
Moreover, $\mu^*$ is the unique $f$-invariant lift of $\nu^*$.
\end{enumerate}
More generally, if $\mu$ is {a hyperbolic non-atomic} ergodic measure {supported on} $\hess$, then $\pi^{-1}(\pi(x))=\{x\}$  $\mu$-a.e.
\end{mainthm}

%Old
%\begin{theoremA}\label{thm:A}
%The essential factor $g$ has a unique MME $\nu^*$, and $\supp(\nu^*)=
%=S_{\mathrm{ess}}$. There is a unique positive-entropy homoclinic class
%$\hess$ of $f$ whose image is not a finite periodic orbit. It satisfies
%$\pi(\hess)=S_{\mathrm{ess}}$ and
%\[
% \htop(g)=\htop(f|_{\hess})
% =\sup\{h_\mu(f):\mu\text{ ergodic and }\mu\hrel\hess\}.
%\]
%Let $\mu^*$ be the MME of $f|_{\hess}$. Then $\nu^*=\pi_*\mu^*$,
%and $\pi$ induces an isomorphism between $(f,\mu^*)$ and $(g,\nu^*)$.
%The measure $\mu^*$ is the unique invariant lift of $\nu^*$.
%Both systems are Bernoulli.
%
%More generally, if $\mu$ is ergodic, $h_\mu(f)>0$, and
%$\mu\hrel\hess$, then $\pi^{-1}(\pi(x))=\{x\}$ for $\mu$-almost
%every $x$. Every other positive-entropy homoclinic class is contained
%in the preimage of a finite periodic orbit of $g$.
%\end{theoremA}

We call $\hess$ the \emph{essential homoclinic class}, relative to the
quotient $\pi$. Essentiality here means that the class is not
collapsed to a periodic orbit. The theorem also shows that its image
is the whole essential factor.

Theorem \ref{thm:A} is consequence of the following statement, which has no
hypothesis on the isotopy class of $f$. Below, $Y$ is a compact metric space.

\begin{mainthm}[\bf Transitive monotone factors]\label{thm:general}
Let $f:S\to S$ be a $C^\infty$ surface diffeomorphism, let $g:Y\to Y$
be a transitive homeomorphism with $\htop(g)>0$, and
let $\pi:S\to Y$ be a monotone semiconjugacy. Then:
\begin{enumerate}[\rm (1)]
\item $f$ has a unique homoclinic class $H$ with positive entropy {and}
such that $\pi(H)$ is infinite. 
%Every other homoclinic class with positive entropy is mapped {\pur under $\pi$} to a finite periodic orbit.
\item If $\mu$ is {a hyperbolic} ergodic measure {supported on} $H$, then $\pi^{-1}(\pi(x))=\{x\}$ $\mu$-a.e. In particular
$\pi:(f,\mu)\to(g,\pi_*\mu)$ is an isomorphism.
\item $\htop(g)=\htop(f|_H)$, and if $\mu^*$ is the MME for $f|_H$ then $\nu^*=\pi_*\mu^*$ is
the unique MME for $g$; both measures are Bernoulli up to the period of $H$,
and $\mu^*$ is the unique $f$-invariant lift of $\nu^*$.
\item If the union of the supports of the {$g$-invariant} ergodic    
measures with positive entropy is dense in $Y$, then $\pi(H)=Y$
and $\nu^*$ has full support.
\item If $g$ is totally transitive, then the period of $H$ is one and
$\nu^*$ is Bernoulli. In particular, this holds when $g$ is mixing.
\end{enumerate}
\end{mainthm}

A map $g$ is totally transitive if $g^n$ is transitive for every $n\neq 0$. We do not know if {the assumption of} item (4) is always verified and, in particular, if $\pi(H)=Y$.
Note that a positive answer ---applied with $Y=S$ and $\pi=\mathrm{Id}$--- would imply that if $f$ is transitive then $H=S$.

\subsection{Origins and related results}
The idea of identifying an essential homoclinic class comes from
Addas-Zanata's work \cite{addaszanata15area}. In the setting
of area-preserving torus diffeomorphisms with rotation sets of nonempty
interior, he constructed hyperbolic periodic saddles whose invariant
manifolds connect integer translates and organize a global region of
instability. This suggests looking for a distinguished homoclinic
class and separating it from the dynamics confined to bounded regions.
Our argument implements this idea through the quotient: it proves
that exactly one class with positive entropy survives, and that all the
others lie in finite unions of periodic fibres. The proof does not
require the extension $f$ to preserve area.

The proof of Theorem \ref{thm:general} uses $su$-quadrilaterals and the
spectral decomposition of \cite{buzzicrovisiersarig22measures}.
Indeed, the existence of the MME inside each
homoclinic class with positive entropy is already proved in their Theorem 2. Buzzi constructed $C^r$ counterexamples to this statement for every $r < \infty$ \cite{buzzi14nomme}. 
Sarig proved that ergodic equilibrium states of H\"older continuous potentials with positive entropy
are Bernoulli times a cyclic permutation \cite{sarig11bernoulli}.

A related recent result of Carrasco, Rodriguez Hertz, Rodriguez Hertz \&
Ures concerns
a smooth diffeomorphism in a pseudo-Anosov isotopy class whose
entropy equals that of the pseudo-Anosov representative \cite{carrascorodriguezhertzrodriguezhertzures26entropy}.
They prove
that Handel's semiconjugacy (which is in this case global by \cite[Theorem 0.1]{handel88entropy}) has cellular fibres and that the unique
lift of the pseudo-Anosov MME is Bernoulli, with
the semiconjugacy giving a measure-theoretic isomorphism. {In the work in preparation \cite{decarvalhogarcia26handel},
the authors prove} that the fibres of the essential factor are the connected components of the fibres of Handel's semiconjugacy, which gives that for the same-entropy case, both semiconjugacies coincide and satisfy Theorem \ref{thm:A}.

In the setting of Theorem \ref{thm:A}, 
monotonicity is a direct consequence of the construction of the essential factor. 
The diffeomorphism $f$ does not need to be isotopic to a genuine pseudo-Anosov map, and equality of
the topological entropy of $f$ and its factor $g$ is not needed.
Theorem \ref{thm:general} isolates the consequences of monotonicity and transitivity
once such a quotient is available.

The distinction between the entropy of a class and that of the
whole extension is also relevant for periodic fibres. 
A homoclinic class contained in such a fibre can have positive entropy {and still be projected to a} periodic orbit. 
Our theorem identifies exactly which
class contributes to the entropy of the transitive quotient.

\subsection{Outline of the argument}
The $su$-quadrilaterals associated to periodic orbits, of arbitrarily
small diameter, cover
almost every point of a hyperbolic measure with positive entropy
(Section~\ref{sec:prelim}). A connected fibre leaving
one of these quadrilaterals meets its boundary, so its image is
forward or backward asymptotic to a periodic orbit. Recurrence
excludes this for a non-atomic projected measure. Shrinking the
quadrilaterals then proves that almost every fibre is a singleton
(Section~\ref{sec:fibres}).

A transitive factor admits a compact transitive lift. Since a
non-atomic $g$-invariant measure with positive entropy has a unique $f$-invariant lift, all
such lifted measures are supported on this same transitive set.
The homoclinic class theory of \cite{buzzicrovisiersarig22measures}
then places them in a single class (Section~\ref{sec:transitivity}).
The unique MME of that class projects without losing entropy and is the
unique MME of the factor. Total transitivity eliminates periods, and
local positive entropy gives full support in the essential factor
(Section \ref{sec:application}).
Section~\ref{sec:consequences} records a principal compact extension
and the resulting existence of equilibrium states for continuous
potentials.

{\subsection{AI use disclaimer} During the production of this text, Garcia-Sassi and Tal used an LLM to assist in checking proofs and correcting the redaction. While the main ideas were all developed without the use of generative AI, some of the suggestions of the model were incorporated in streamlining proofs and, in case of Theorem B, simplifying hypotheses.}

\section*{Acknowledgements}

The second author was supported by FAPESP grant 2025/08777-1, IMERL and PEDECIBA. The third author was supported by CNPq/Universal 401737/2025-0, Instituto Serrapilheira grant “Jangada Dinâmica: Impulsionando Sistemas Dinâmicos na Região Nordeste”, and FAPESP grant 2025/11400-7. The fourth author was supported by CNPq and FAPESP grant 2023/07076-4.

\section{Hyperbolic measures and $su$-quadrilaterals}\label{sec:prelim}

For a saddle periodic orbit
$\mathcal O$, write $W^{s/u}(\mathcal O)=\bigcup_{p\in\mathcal O}W^{s/u}(p)$ {for its stable and unstable manifolds}.
Two saddle periodic orbits $\mathcal O,\mathcal O'$ are \emph{homoclinically related} if each {set of transverse intersections}
$W^s(\mathcal O)\pitchfork W^u(\mathcal O')$ and
$W^u(\mathcal O)\pitchfork W^s(\mathcal O')$ {are non-empty}. When this happens,
we write $\mathcal O\hrel\mathcal O'$.
The homoclinic class of $\mathcal O$ is {the transitive invariant compact set} defined as
\[
 \HC(\mathcal O)=\overline{\bigcup\{\mathcal O':\mathcal O'\hrel\mathcal O\}}=\overline{W^s(\mathcal O)\pitchfork W^u(\mathcal O)}.
\]
The second equality follows from the inclination lemma.

\begin{definition}[\bf Homoclinic relation of measures \cite{buzzicrovisiersarig22measures}]\label{def:homoc_relat_meas}
	Given two ergodic hyperbolic measures $\mu,\nu$,
	we say that {\em $\mu$ and $\nu$ are homoclinically related} if for $\mu$-a.e. $x$ and
	$\nu$-a.e. $y$ it holds 
	$$W^s(\mathcal O_x)\pitchfork W^u(\mathcal O_y)\neq \emptyset \ \text{ and }\ 
	{W^u(\mathcal O_x)\pitchfork W^s(\mathcal O_y)}\neq \emptyset,
	$$
	where $W^{s/u}(\mathcal O_x)$ denotes the global stable manifold of the orbit of $x$,
	and similarly for $y$.
	
\end{definition}

\medskip
Given a {saddle} period orbit $\mathcal O$, we will write
$\mu\hrel\mathcal O$ when $\mu$ is homoclinically related to the unique $f$-invariant
probability measure supported on $\mathcal O$.
The homoclinic relation is an equivalence relation \cite[Prop. 2.11]{buzzicrovisiersarig22measures}.

Recall that the period of $\HC(\mathcal{O})$ is
$\ell(\mathcal O)=\gcd\{|\mathcal O'|:\mathcal O'\text{ { saddle periodic orbit,} }\mathcal O'\hrel\mathcal O\}$.
The following theorem summarizes properties we will need.

\begin{theorem}[\cite{buzzicrovisiersarig22measures}]\label{prop:BCS}
If $f$ is a $C^\infty$ surface diffeomorphism, then the following holds.
\begin{enumerate}[{\rm (1)}]
\item Every {hyperbolic, non-atomic,} ergodic probability measure $\mu$ is homoclinically related to a saddle periodic orbit $\mathcal O$, and $\supp(\mu)\subset\HC(\mathcal O)$. {Conversely,} if $\mu$ is {hyperbolic, non-atomic,} ergodic, and supported on ${\HC(\mathcal O)}$,
then it is homoclinically related to $\mathcal O$.
\item Every homoclinic class $H=\HC(\mathcal O)$ {supports} a unique
{measure} $\mu_H$ such that
$h_{\mu_H}(f)=h_{\rm top}(f|_{H})$. Moreover, $\supp(\mu_H)=H$ and $\mu_H$ is Bernoulli times
a cyclic permutation of order $\ell(\mathcal O)$.
\item Every homoclinic class $H=\HC(\mathcal O)$ decomposes as a union $H=H_1\cup\dots\cup H_\ell$ 
with $\ell:=\ell(\mathcal O)$,
such that $H_1,\dots,H_\ell$ are distinct homoclinic classes of $f^\ell$ which are cyclically permuted by $f$.
\item If $K$ is a compact transitive $f$-invariant set, then every two ergodic positive entropy
measures supported on $K$ are homoclinically related.
\end{enumerate}
\end{theorem}

Part (1) follows from the work of Katok \cite{katok80lyapunov},
see also \cite[Cor. 6.4.19]{katokhasselblatt95introduction}, and it holds when $f$ is just $C^{1+\alpha}$.
The second part is \cite[Theorem 2(1)]{buzzicrovisiersarig22measures},
{the third part is \cite[Theorem 6.12(3)]{buzzicrovisiersarig22measures},}
and the fourth part is \cite[Theorem 1(5)]{buzzicrovisiersarig22measures}.
The proof of part (4) indeed only uses the existence of $su$-quadrilaterals for measures
with positive entropy and the dynamical Sard theorem \cite[Theorem 4.2]{buzzicrovisiersarig22measures}.

\medskip
\noindent {\sc $su$-quadrilateral:}
An {\em $su$-quadrilateral} associated to a saddle periodic orbit $\mathcal O$ is
an open topological disk $Q$ whose boundary is a simple closed curve
formed by four compact arcs, alternately contained in $W^s(\mathcal O)$ and
$W^u(\mathcal O)$, with transverse corners.

\begin{lemma}\label{lem:quads}
Let $f$ be a $C^{1+\alpha}$ {surface diffeomorphism}, let $\mu$ be {a hyperbolic, non-atomic,} ergodic measure, and let $\mathcal O$ be a {saddle} periodic orbit
homoclinically related to $\mu$.
\begin{enumerate}[{\rm (1)}] 
\item For all $\varepsilon>0$, there is a finite
family of $su$-quadrilaterals $Q_1,\ldots,Q_N$ associated to $\mathcal O$
such that
\[
 \diam(Q_i)<\varepsilon\text{ for }i=1,\ldots,N,\text{ and }
 \mu\left(\bigcup_{i=1}^NQ_i\right)>1-\varepsilon.
\]
\item There is a countable family $\mathcal Q$ of $su$-quadrilaterals associated to $\mathcal O$ such that
for $\mu$-a.e. $x$ there are $Q_1,Q_2,\ldots\subset \mathcal Q$ such that
$x\in Q_i$ for all $i\geq 1$ and $\diam(Q_i)\to 0$ as $i\to+\infty$.
\end{enumerate}
\end{lemma}

\begin{proof}
(1) Proposition 2.19 of \cite{buzzicrovisiersarig22measures} establishes a more general result,
allowing $\mu$ to be non-ergodic and assuming that $\mu$ has no atoms.
First take $su$-quadrilaterals $\widetilde Q_1,\dots,\widetilde Q_N$ associated to periodic
orbits $\mathcal O_1,\dots,\mathcal O_N$ homoclinically related to hyperbolic ergodic components of $\mu$.
When $\mu$ is ergodic, all the periodic orbits $\mathcal O_i$
are homoclinically related. One fixes one of them, say $\mathcal O$, and the inclination
lemma shows that any $su$-quadrilateral $\widetilde Q_i$ can be replaced by some
$su$-quadrilateral $Q_i$ associated to $\mathcal O$ such that $\diam(Q_i)$ and
$\diam(\widetilde Q_i)$ are close and $\mu(Q_i \triangle \widetilde Q_i)$ is small. Hence item (1) follows.\\

%We explain why the orbit can be prescribed. The inclination lemma
%\cite[Lemma 2.7]{buzzicrovisiersarig22measures} implies that each
%compact stable arc in $W^s(O_i)$ can be approximated in $C^1$ by
%arcs in $W^s(O)$; the same holds for unstable arcs. Approximate
%the four sides of each quadrilateral simultaneously, first extending
%each side slightly past its corners. The transverse
%corners persist, and the selected intersections of adjacent sides
%bound a nearby quadrilateral associated to $O$. This is the same
%replacement used in \cite[proof of Theorem 6.4, Step 1]
%{buzzicrovisiersarig22measures}.

%For the measure estimate, first choose the original diameters
%less than $\varepsilon/2$ and the measure of their union greater
%than $1-\varepsilon/2$. Choose compact subsets of their interiors
%whose union has measure greater than $1-\varepsilon$. Sufficiently
%close replacements contain these compact subsets and have diameters
%less than $\varepsilon$. This gives the asserted family.

\noindent
(2) Applying part (1) for $\varepsilon_j=2^{-j}$, we obtain a 
finite family $\mathcal Q_j$ of $su$-quadrliaterals whose complement has measure $<2^{-j}$.
Define $\mathcal Q=\bigcup\mathcal Q_j$, which is a countable family of $su$-rectangles whose
union has full $\mu$-measure. 
Since $\sum\mu(\mathcal Q_j^c)<\infty$, the Borel-Cantelli lemma 
implies that for $\mu$-a.e. $x$ there is $j_0=j_0(x)$ such that $x\in \mathcal Q_j$ for all
$j\geq j_0$. This completes the proof.
\end{proof}

\section{A dichotomy for homoclinic classes}\label{sec:fibres}

In this section we let $f:S\to S$ be a $C^{1+\alpha}$ {surface} diffeomorphism, $g:Y\to Y$ be a
homeomorphism on a compact metric space $Y$, and $\pi:S\to Y$ be a {\em monotone
semiconjugacy}, meaning that it is a continuous surjective map with $\pi\circ f=g\circ \pi$
such that $\pi^{-1}(y)$ is connected for every $y\in Y$.
No transitivity or isotopy hypothesis is imposed.

We let $\mu$ be a $f$-invariant ergodic probability measure. The projection
$\nu=\pi_*\mu$, which is a $g$-invariant ergodic probability measure,
is either non-atomic or it is a periodic measure (supported on a periodic orbit for $g$). In the sequel,
we characterize the projection of homoclinic classes
under this dichotomy.

\begin{theorem}\label{thm:singleton}
Let $\mu$ be a {hyperbolic non-atomic} $f$-invariant ergodic probability measure.
\begin{enumerate}[{\rm (1)}]
\item If $\nu=\pi_*\mu$ is non-atomic, then
$\pi^{-1}(\pi(x))=\{x\}$ $\mu$-a.e., $\pi$ induces an isomorphism between $(f,\mu)$ and $(g,\nu)$,
and $\mu$ is the unique Borel probability
on $S$ that projects to $\nu$. In particular, $h_\mu(f)=h_\nu(g)$.
\item If $\nu=\pi_*\mu$ is supported on a periodic orbit $\mathcal O'$,
then for every periodic orbit $\mathcal O$ homoclinically related to $\mu$
it holds $\pi(\HC(\mathcal O))=\mathcal O'$.
\end{enumerate}
\end{theorem}

\begin{proof}
(1) Choose a periodic orbit $\mathcal O$ homoclinically related to $\mu$.
Let $\mathcal Q$ be a countable family of $su$-quadrilaterals associated to $\mathcal O$
satisfying Lemma \ref{lem:quads}, and let $X_0$ be the full $\mu$-measure set of points
$x\in S$ for which there are $Q_1,Q_2,\ldots\subset\mathcal Q$ such that
$x\in Q_i$ and $\diam(Q_i)\to 0$.
Define
$$
Y_0=\{y\in Y: \exists m_k,n_k\to+\infty \text{ s.t. }
g^{m_k}(y),g^{-n_k}(y)\to y\text{ as }k\to+\infty\},
$$
which has full $\nu$-measure by the Poincaré recurrence theorem.
We let $\mathcal O'=\pi(\mathcal O)$, which has zero $\nu$-measure.
We claim that $\pi^{-1}(y)$ is a singleton for $\nu$-a.e. $y\in Y_0$. If this is not the case,
there is $y\in Y_0$, $y\notin\mathcal O'$, such that $\# \pi^{-1}(y)>1$ and
$\pi^{-1}(y)\cap X_0\neq\emptyset$. Let $x\in \pi^{-1}(y)\cap X_0$ and $\widetilde x\in \pi^{-1}(y)\setminus\{x\}$.
Let $Q_1,Q_2,\ldots\subset\mathcal Q$ such that
$x\in Q_i$ and $\diam(Q_i)\to 0$.
If $i$ is large, then $\widetilde x\notin \overline{Q_i}$. Since $\pi^{-1}(y)$ is connected and 
intersects $Q_i$ and $\overline{Q_i}^c$, it necessarily intersects $\partial Q_i$. 
Assuming without loss of generality that such intersection contains an element
$z\in W^s(\mathcal O)$, we conclude that $y=\pi(z)$ is forward asymptotic to $\mathcal O'$,
hence it is not forward recurrent. {This gives a contradiction.}

The rest of the proof is standard. The set 
\[
Y_1=\{y\in Y:\pi^{-1}(y)\text{ is a singleton}\}
 =\bigcap_{n\ge1}\{y\in Y:\diam(\pi^{-1}(y))<1/n\}
\]
is measurable, hence $X_1=\pi^{-1}(Y_1)$ is also measurable. Since
$\pi|_{X_1}:X_1\to Y_1$ is a measurable bijection, the Lusin-Suslin
theorem implies that its inverse is measurable. 
We proved that $\nu(Y_1)=1$, hence $\pi|_{X_1}$ is a measurable
isomorphism between $(f,\mu)$ and $(g,\nu)$. In particular,
$h_\mu(f)=h_\nu(g)$. Finally, every lift of $\nu$ to $S$
is carried by $X_1$, hence $\mu$ is the unique lift.\\

\noindent
(2) Since $HC(\mathcal O)$ is the closure of the set of periodic orbits $\widetilde{\mathcal O}$
homoclinically related to $\mathcal O$, we just need to prove that 
$\pi(\widetilde{\mathcal O})=\mathcal O'$ for every such orbit. 
By contradiction, assume that $\widetilde{\mathcal O}$ is homoclinically
related to $\mathcal O$ (hence also to $\mu$) and satisfies $\pi(\widetilde{\mathcal O})\ne\mathcal O'$.
We can fix $\widetilde{\mathcal Q}$ a countable
family of $su$-quadrilaterals associated to $\widetilde{\mathcal O}$ satisfying  
Lemma \ref{lem:quads}, and apply the same argument used in part (1)
to conclude that 
$\pi^{-1}(y)$ is a singleton for every $y\in\mathcal O'$, otherwise
we would find $z\in \pi^{-1}(y)\cap W^{s/u}(\widetilde{\mathcal O})$ and
so $y=\pi(z)$ would not be in $\mathcal O'$. But then $\mu$ is supported
on the finite set $\pi^{-1}(\mathcal O')$, {which contradicts the assumption that
$\mu$ is non-atomic}. This concludes the proof.
\end{proof}

\begin{corollary}\label{cor:dichotomy}
Let {$f:S\to S$} be a $C^{1+\alpha}$ surface diffeomorphism, $g:Y\to Y$ a
homeomorphism on a compact metric space, and $\pi:S\to Y$ a monotone
semiconjugacy. For each homoclinic class $H=\HC(\mathcal{O})$, exactly one of the following holds:
\begin{enumerate}[{\rm (1)}]
\item $\pi(H)$ is a finite periodic orbit $\mathcal O'$, and every $f$-invariant
measure supported on $H$ projects to the periodic measure on $\mathcal O'$;
\item $\pi(H)$ is infinite, and every {hyperbolic non-atomic} ergodic
measure $\mu$ supported on $H$ projects to a
$g$-invariant ergodic measure $\nu$ with same entropy and 
singleton fibres $\nu$-a.e.
\end{enumerate}
\end{corollary}

\begin{proof}
(1) Assume that $\pi(H)$ is finite. Since $f|_H$ is transitive, the same happens
to $g|_{\pi(H)}$, and so $\pi(H)$ is a periodic orbit.

\medskip
\noindent
(2) Assume that $\pi(H)$ is infinite. Let $\mu$ be {hyperbolic non-atomic} ergodic with $\mu(H)=1$. By Theorem~\ref{prop:BCS}(1), {the periodic orbit $\mathcal O$ is}
homoclinically related to 
$\mu$. By Theorem \ref{thm:singleton}(2), $\pi_*\mu$ is non-atomic, otherwise
$\pi(H)$ would be finite. Now apply Theorem~\ref{thm:singleton}(1).
\end{proof}

\section{Proof of Theorem \ref{thm:general}}\label{sec:transitivity}

We begin with some general facts, that do not require $f$ to be $C^\infty$.

%\begin{lemma}[Invariant lifts]\label{lem:lifting}
%Let $\pi:(X,f)\to(Y,g)$ be a factor map between homeomorphisms
%of compact metric spaces. Every ergodic invariant probability
%$\nu$ of $g$ has an ergodic invariant lift $\mu$ with
%$h_\mu(f)\ge h_\nu(g)$.
%\end{lemma}
%
%\begin{proof}
%There is a probability $\lambda$ on $X$ with $\pi_*\lambda=\nu$.
%For example, approximate $\nu$ by finite atomic probabilities,
%choose one preimage of each atom, and pass to a convergent
%subsequence of the lifted probabilities. Any weak-$*$ limit of
%$N^{-1}\sum_{j=0}^{N-1}f_*^j\lambda$ is invariant and projects
%to $\nu$. In its ergodic decomposition, almost every component
%projects to $\nu$, since $\nu$ is an extreme point of
%$\mathcal M_g$. Choose one such component. The entropy inequality
%is the usual factor inequality.
%\end{proof}

\begin{lemma}\label{lem:transitive-lift}
{Let $f:S\to S$ be a $C^{1+\alpha}$ surface diffeomorphism, $g:Y\to Y$ a
homeomorphism on a compact metric space,} $\pi:S\to Y$ a semiconjugacy
{and $Y'\subset Y$ a $g$-invariant compact set such that $g|_{Y'}$ is transitive. Then:}
\begin{enumerate}[{\rm (1)}]
\item {There exists a $f$-invariant compact set $K$} {such that}
$\pi(K)=Y'$ {and} such that $f|_{K}$ is transitive. 
\item If $\pi$ is monotone, then for every {hyperbolic} ergodic probability measure
$\mu$ such that $\pi_*\mu$ is non-atomic {and supported on $Y'$}, it holds 
$\supp(\mu)\subset K$.
\end{enumerate}
\end{lemma}

\begin{proof}
(1) This property is general and only requires continuity. Consider the
family
$$
\mathcal Z=\{X'\subset X: X'\text{ is compact $f$-invariant with } \pi(X')\supset Y'\}
$$
partially ordered
by inclusion. By the Zorn lemma, there is $X'$ minimal. Let $y\in Y'$ with dense orbit
on $Y'$, and let $x\in X'$ projecting to $y$. Then $\omega(x)\subset X'$ is a compact $f$-invariant
set such that $\pi(\omega(x))=\omega(y)=Y'$. By minimality,
$\omega(x)=X'$ and $\pi(X')=\pi(\omega(x))=Y'$. 

%Choose a compact invariant $K$ minimal under inclusion among those
%mapping onto $Y$. Such a $K$ exists by the nested-intersection
%argument and Zorn's lemma: the intersection of a chain remains
%invariant, and its intersection with each compact fibre is
%nonempty.

%Birkhoff's transitivity theorem, applied in both time directions,
%provides $y\in Y$ with $\omega_g(y)=\alpha_g(y)=Y$. One can
%see this by intersecting the countably many open dense sets of
%points whose forward, respectively backward, orbit visits each
%basis open set after any prescribed time. Choose $x\in K$ over
%$y$. The compact invariant sets $\omega_f(x)$ and $\alpha_f(x)$
%both map onto $Y$, because continuous factor maps commute with
%alpha- and omega-limit sets on compact spaces. Minimality gives
%$\omega_f(x)=\alpha_f(x)=K$, proving transitivity.

\medskip
\noindent
(2) {Since $\nu=\pi_*\mu$ is non-atomic, $\mu$ is non-atomic and Theorem~\ref{thm:singleton}(1) applies:} $\mu$ is the unique lift of $\nu$.
Since $\pi|_{K}:K\to Y'$ {is a semiconjugacy}, $\nu$ has a lift to $K$, and so
$\supp(\mu)\subset K$.
\end{proof}

We are now ready to prove Theorem \ref{thm:general}.

\begin{proof}[Proof of Theorem \ref{thm:general}] (1)
%The connected space $S$ maps onto $Y$, so $Y$ is connected.
%Positive entropy excludes a one-point space; hence $Y$ has no
%isolated points. 
Let $K$ be the set given by Lemma~\ref{lem:transitive-lift}(2).
By the variational principle, there is a $g$-invariant ergodic {measure} $\nu$ with $h_\nu(g)>0$.
Lift $\nu$ to an $f$-invariant ergodic measure $\mu$ on $K$. Hence 
$h_\mu(f)>0$ and so there is a periodic orbit $\mathcal O$ homoclinically related to $\mu$.
Define $H:={\rm HC}(\mathcal O)$.
Since $\supp(\mu)\subset H$ and $\nu=\pi_*\mu$ has positive entropy, 
$\pi(H)$ is infinite.

We now prove that every other homoclinic class projects to a periodic orbit.
We begin proving that every two $f$-invariant ergodic {measures} $\mu,\mu'$ with positive
entropy whose projections $\pi_*\mu,\pi_*\mu'$ are non-atomic
are necessarily homoclinically related. Indeed, Lemma \ref{lem:transitive-lift}(2)
implies that the supports of $\mu,\mu'$ are {included} in $K$. By Theorem~\ref{prop:BCS}(4),
$\mu,\mu'$ are homoclinically related. 

Now, let $H'$ be a homoclinic class such that
$\pi(H')$ is infinite. By Theorem~\ref{prop:BCS}(2), there is $\mu_{H'}$ 
a measure of maximal entropy for $f|_{H'}$. This measure has support $H'$,
hence its projection has support $\pi(H')$, hence is non-atomic. Applying this
for two possible homoclinic classes with infinite projection, the claim in the previous
paragraph implies that {their measures maximizing the entropy} are homoclinically related, and hence the classes
coincide. This proves that $H$ is the unique homoclinic class with infinite projection.

\medskip
\noindent
(2) Follows directly from Corollary \ref{cor:dichotomy}(2).

\medskip
\noindent
(3) If $\nu$ is $g$-invariant with $h_\nu(g)>0$, then $\nu(\pi(H))=1$
and so $\nu$ lifts to $\mu$ such that $\mu(H)=1$ and $h_\mu(f)\geq h_\nu(g)$.
Then $h_{\rm top}(g)\leq h_{\rm top}(f|_H)$. By Theorem~\ref{prop:BCS}(2),
$f|_H$ has a unique MME $\mu^*$. By Corollary \ref{cor:dichotomy}(2), $\pi|_H$ conjugates $(f,\mu^*)$
and $(g,\nu^*)$, hence $h_{\rm top}(g)=h_{\rm top}(f|_H)$
and $\nu^*$ is a MME for $g$. If $\nu$ were another MME for $g$, then its lift
to $H$ would give another MME for $f|_H$, contradicting Theorem~\ref{prop:BCS}(2).
Finally, also by Theorem~\ref{prop:BCS}(2), $\nu^*$ is Bernoulli times a cyclic permutation
of the period of $H$.

\medskip
\noindent
(4) If $\nu$ is an ergodic $g$-invariant measure with positive entropy,
then $\nu(\pi(H))=1$ and so it lifts to $H$. Hence 
$\supp(\nu)\subset\pi(H)$. Taking the union over all ergodic $g$-invariant measure with positive entropy, 
we get that $\pi(H)=Y$. Since the support of $\mu^*$ is $H$, the support of $\nu^*$ is $Y$.

\medskip
\noindent
(5) Let $\ell$ be the period of $H$. According to Theorem~\ref{prop:BCS}(3), $H$ decomposes as
$H=H_1\cup\dots\cup H_\ell$. Each class $H_i$ supports a unique measure $\mu_i$
maximizing the entropy of the class for $f^\ell$ and the measure $\mu^*$ decomposes as 
$\mu^*=\frac1\ell\sum_{j=0}^{\ell-1}\mu_j$. Each set $\pi(H_i)$ has positive entropy for $g^\ell$,
hence is infinite. Since $g$ is totally transitive, Theorem~\ref{thm:general}(1) applies {to $f^\ell$}
and so there exists only one class $H_i$. This implies that $\ell=1$, i.e. $\nu^*$ is Bernoulli.
\end{proof}

\section{Proof of Theorem \ref{thm:A}}\label{sec:application}

The \textit{essential factor} constructed in \cite[Theorems A and B]{garciatal26essential} for a surface homeomorphism in a relative pseudo-Anosov isotopy class is obtained by collapsing maximal dynamically bounded continua. 
For that, the authors use a generalized stable-unstable foliation by \cite{militon24foliations} to prove that the lifted diameter of this sets is uniformly bounded. This yields a monotone cellular semiconjugacy, which is homeomorphic to the original surface \cite{moore25continua, robertssteenrod38monotone}. By construction, the factor dynamics $g$ is infinitely continuum-wise expansive and has the following properties:
\begin{enumerate}[$\circ$]
\item $g$ is topologically mixing, hence it is totally transitive.
\item Every nonempty open set $U\subset Y$ contains a topological horseshoe, which gives a
compact set $\Lambda$ invariant under some $g^n$ with
$\htop(g^n|_\Lambda)>0$. If $\nu$ is a $g^n$-invariant measure of positive 
entropy with $\nu(\Lambda)=1$, then $\tfrac{1}{n}\sum_{j=0}^{n-1}g^j_*\nu$ is a $g$-invariant measure
of positive entropy whose support intersects $U$. Hence item (4) in Theorem
\ref{thm:general} holds.
\end{enumerate}
Theorem \ref{thm:general} then directly implies Theorem \ref{thm:A},
with $\hess=H$.

\begin{remark}[Entropy removed by the quotient]
In general the appropriate equality is
$\htop(g)=\htop(f|_{\hess})\le\htop(f)$, since
a class with positive entropy that is mapped to a periodic orbit makes no contribution to the
entropy of $g$. We thus do not require $\htop(g)=\htop(f)$
as assumed in \cite{carrascorodriguezhertzrodriguezhertzures26entropy}.
\end{remark}

\section{A principal factor and equilibrium states}\label{sec:consequences}

A factor map is called \emph{principal} if it preserves the metric
entropy of every invariant probability. The factor
$\pi:S\to Y$ needs not be principal, since it may collapse
classes with positive entropy. The compact transitive lift used in the
proof {of Theorem \ref{thm:general}} provides a principal restriction.
{Write $\mathcal M_f$ for the $f$-invariant Borel probabilities (same for $g$).}

\begin{corollary}\label{cor:principal}
Under the {setting} of Theorem \ref{thm:general}, {let $K$ be a transitive $f$-invariant compact set such that $\pi(K)=Y$,
cf. Lemma~\ref{lem:transitive-lift}.} Then:
\begin{enumerate}[{\rm (1)}]
\item $K$ contains $H$, and $\pi:(K,f|_K)\to(Y,g)$ is principal. 
\item The entropy map
$\nu\in\mathcal M_g \mapsto h_\nu(g)$ is upper semicontinuous, hence
every continuous potential $\psi:Y\to\R$ has an
equilibrium state.
\end{enumerate}
\end{corollary}

{We remark that, under the above assumptions, equilibrium states might not
be unique.}

\begin{proof}
(1) The measure $\mu^*$ {maximizing the entropy on $H$} has non-atomic projection, so
Lemma~\ref{lem:transitive-lift} gives that $\supp(\mu^*)=H\subset K$.
Let $\mu$ be an ergodic $f$-invariant measure supported on $K$.
If $h_\mu(f)>0$, Theorem~\ref{prop:BCS}(4) applied to $\mu$
and $\mu^*$ on the compact transitive set $K$ gives that
$\mu\hrel H$, and {so Theorem \ref{thm:general}(2) implies that
$h_\mu(f)=h_{\pi_*\mu}(g)$. If $h_\mu(f)=0$ then the same equality holds.}
Integration over ergodic components proves that
$h_\mu(f)=h_{\pi_*\mu}(g)$ for every $\mu\in\mathcal M_f$ with $\mu(K)=1$.\\

\noindent
(2) Let $\{\nu_n\} \subset\mathcal M_g$ converging to $\nu$ in the weak-* topology.
For each $n$, let $\mu_n$ be a lift of $\nu_n$ supported on $K$. Passing to two subsequences,
we find a subsequence $\{n_k\}$ such that
$\limsup\limits_{n\to+\infty} h_{\nu_n}(g)=\lim\limits_{k\to+\infty} h_{\nu_{n_k}}(g)$ and
$\mu_{n_k}$ converges to $\mu\in\mathcal M_f$ in the
weak-* topology. Then $\supp(\mu)\subset K$ and $\pi_*\mu=\nu$.
Since $f$ is $C^\infty$, the entropy map of $f$ is upper semicontinuous
\cite{newhouse89continuity}; see also
\cite[Section 5.1]{buzzicrovisiersarig22measures}. Hence
\[
 \limsup_{n\to+\infty} h_{\nu_n}(g)=\lim\limits_{k\to+\infty} h_{\nu_{n_k}}(g)=\limsup_{k\to+\infty} h_{\mu_{n_k}}(f)
 \le h_\mu(f)=h_\nu(g).
\]
Therefore, the function
\[
 \nu\in\mathcal M_g\longmapsto h_\nu(g)+\int_Y\psi\,d\nu
\]
is upper semicontinuous and so it attains a maximum.
\end{proof}

{We say a map is \emph{almost one-to-one} if it is injective on a dense subset of its domain.}

\begin{corollary}\label{cor:minimal-core}
Under the setting of Theorem~\ref{thm:general}(4),
\[
 Y_1=\{y\in Y:\pi^{-1}(y)\text{ is a singleton}\} \ \text{ and } \ X_1=\pi^{-1}(Y_1)
\]
are invariant dense $G_\delta$ subsets of $Y$ and $H$, respectively.
In particular, {$\pi:(H,f|_H)\to(Y,g)$} is a principal almost one-to-one factor.
\end{corollary}

\begin{proof}
{By Theorem~\ref{thm:general}(4), $\pi(H)=Y$. Applying Corollary~\ref{cor:principal} to $K=H$,
we get that ${\pi:(H,f|_H)\to(Y,g)}$ is principal.

The set $Y_1$ is $G_\delta$ by the upper semicontinuity argument in Theorem~\ref{thm:singleton}.
The measure $\mu^*$ maximizing the entropy on $H$ has full support, hence the support of $\nu^*=\pi_*\mu^*$ is $Y$ and by Theorem~\ref{thm:singleton} one concludes that $Y_1$ is dense in $Y$.}

Since $\pi(H)=Y$, every singleton fibre belongs to $H$. Hence $X_1\subset H$ is an
invariant $G_\delta$ with full $\mu^*$-measure, and is dense in
$H$. The inverse of $\pi|_{X_1}$ is continuous: if $y_n\to y$
in $Y_1$,  continuity makes every limit point of their preimages
belong to the unique fibre over $y$. Thus the preimages converge.
\end{proof}

%\begin{remark}
%Corollary~\ref{cor:principal} asserts existence for all continuous
%potentials, without asserting uniqueness. Uniqueness for
%nonconstant potentials requires additional control. In particular,
%a continuous or H\"older potential on $Y$ need not have a
%H\"older pullback under an arbitrary monotone map $\pi$.
%\end{remark}

\bibliographystyle{plain-like-initial}
\bibliography{biblio}
\end{document}